\documentclass[12 pt]{amsart}
\usepackage[margin=1in]{geometry}
\usepackage{amssymb, amsmath, amsthm, amsfonts}
\usepackage{hyperref}
\usepackage{setspace}
\usepackage[all]{xy}
\usepackage{graphicx}

\newtheorem{theorem}{Theorem}[section]

\newtheorem{proposition}[theorem]{Proposition}

\theoremstyle{definition}
\newtheorem{remark}[theorem]{Remark}
\newtheorem{definition}[theorem]{Definition}

\newcommand{\pullbackcorner}[1][dr]{\save*!/#1-1.7pc/#1:(-1,1)@^{|-}\restore}

\def\beq{\begin{eqnarray*}}
\def\eeq{\end{eqnarray*}}

\def\R{\mathbb{R}}
\def\Z{\mathbb{Z}}

\def\D{\Delta}
\def\incl{\hookrightarrow}
\def\to{\rightarrow}

\def\eps{\varepsilon}

\def\dim{\mathrm{dim}\>}

\def\x{\times}
\def\d{\partial}
\def\phi{\varphi}

\def\Emb{\mathrm{Emb}}

\def\P{\mathcal{P}}

\def\LL{\mathcal{L}}
\def\K{\mathcal{K}}

\def\E{E^\nu}

\def\Bl{\mathrm{Bl}}

\def\Map{\mathrm{Map}}

\def\V{\mathcal{V}}
\def\CD{\mathcal{CD}}
\def\TD{\mathcal{TD}}
\def\bx{\mathbf{x}}

\def\bt{\mathbf{t}}
\def\la{\langle}
\def\ra{\rangle}

    \title[Homotopy Bott--Taubes and the Taylor tower for spaces of knots and links]{Homotopy Bott--Taubes integrals and the Taylor tower for spaces of knots and links}
    \author{Robin Koytcheff}
    \thanks{This work was supported  partially by NSF grant DMS-1004610 and partially by PIMS}
    \address{Mathematics and Statistics, University of Victoria, Victoria, BC, Canada}
    \email{rmjk@uvic.ca}

\begin{document}
\maketitle

\begin{abstract}
This work continues the study of a homotopy-theoretic construction of the author inspired by the Bott--Taubes integrals.  Bott and Taubes constructed knot invariants by integrating differential forms along the fiber of a bundle over the space of knots in $\R^3$.  Their techniques were later used by Cattaneo et al.~to construct real ``Vassiliev-type" cohomology classes in the space of knots in $\R^d$, $d\geq 4$.  By doing this integration via a Pontrjagin--Thom construction, we constructed cohomology classes in the knot space with \emph{arbitrary} coefficients.  
We later showed that a refinement of this construction recovers the Milnor triple linking number for string links.  We conjecture that we can produce all Vassiliev-type classes in this manner.
Here we extend our homotopy-theoretic constructions to the stages of the Taylor tower for the embedding space, which arises from the Goodwillie--Weiss embedding calculus.  
We use the model of ``punctured knots and links" for the Taylor tower.
\end{abstract}

\section{Introduction}
This paper concerns spaces of knots and links.  More precisely, for $d\geq 3$, let $\LL_m^d:=\Emb(\coprod_{i=1}^m \R, \R^d)$, the space of $m$-component long links (a.k.a.~string links) in $\R^d$, that is, the space of embeddings of $m$ disjoint copies of $\R$ into $\R^d$ with prescribed behavior outside of the intervals $\coprod_{i=1}^m [-1,1]\subset \coprod_{i=1}^m \R$.  
This is an obvious generalization of the space $\LL_1^d=\Emb(\R, \R^d)$ of long knots in $\R^d$, which we often denote $\K^d$ or just $\K$.
For $d=3$ and $m=1$, the connected components of $\LL_1^3$ correspond to isotopy classes of long knots, which correspond (bijectively) to isotopy classes of closed knots.
When $d\geq 4$, $\Emb(\coprod_m \R, \R^d)$ is connected but is far from being topologically trivial. (For example, see \cite{Cattaneo, LTV} for the case $m=1$.)  

We build on our previous work \cite{Rbo}, in which we produced cohomology classes with arbitrary coefficients in these spaces via a homotopy-theoretic construction inspired by the configuration space integrals of Bott and Taubes.  In \cite{HoMilnor3ple}, we showed that a ``gluing refinement" of this construction recovers the triple linking number $\mu_{123}$ for string links.  We conjecture that we can produce integral multiples of all the Vassiliev-type classes of Cattaneo--Cotta-Ramusino--Longoni \cite{Cattaneo}, thus showing that these classes are rational.  This conjectured statement is work in progress.
(Our idea there is that graph cocycles come from cancellations of terms in the graph complex;  these cancellations can be grouped into pairs; and these canceling pairs correspond to gluing pairs of configuration spaces, a construction which should be amenable to singular (co)homology and not just de Rham theory.)

The core idea of this paper is to extend the methods of Voli\'c from \cite{VolicConfig} to relate our homotopy-theoretic integration to the embedding calculus of Goodwillie and Weiss.  The embedding calculus produces a Taylor tower of spaces 
\[
...\to T_n \K^d \to T_{n-1} \K^d \to ... \to T_0 \K^d =*
\]
together with compatible maps $\K^d \to T_n \K^d$.  For $d\geq 4$, its inverse limit is the knot space $\K^d$ itself, as shown in work of Goodwillie, Klein, and Weiss \cite{GW, GK}.  Voli\'c's main result in \cite{VolicConfig} is that for $d=3$, all finite-type invariants factor through the tower.  For the space $\LL_m^d$ of $m$-component links in $\R^d$, there is a similar \emph{multi-tower} of spaces, indexed by $m$ nonnegative integers $(n_1,...,n_m)$ rather than one nonnegative integer $n$.  For $d\geq 4$, the work of Goodwillie, Klein, and Weiss also implies convergence of this tower to the link space.  

Our first main result is Theorem \ref{FTLinkInvtsFactor}, which is a generalization from knots to links of Voli\'{c}'s result that finite-type invariants factor through the tower.
The other main result of this paper is that both the basic homotopy-theoretic construction of \cite{Rbo}, as well as the refined construction which yields $\mu_{123}$, extend to the Taylor tower.  These results appear as Proposition \ref{RboOnTaylor} and Theorem \ref{HoMilnorOnTaylor}.  
While the core idea of this paper applies to any $d$ and $m$, we will at some point set $d=3, m=3$ to address the construction of $\mu_{123}$.
We use a ``punctured links" model for the Taylor tower that has appeared in work of other authors, e.g.\cite{Weiss, SinhaTop, VolicConfig}.
Ultimately, we deduce that $\mu_{123}$ factors through the stage $T_{(2,2,2)} \LL_3^3$ of the multi-tower.  

We could deduce this latter result via integration of forms, using 
Theorem \ref{FTLinkInvtsFactor}, which produces all finite-type link invariants in the tower by integration.  The advantage of using a more homotopy-theoretic approach is that one may be able to construct integer or mod $p$ classes that cannot be realized by integration.  Then similar arguments to the ones in this paper would allow one to see those classes in the tower.  

Finally, we point out that the triple linking number for closed links is seen at stage $(1,1,1)$ in the tower for the space of \emph{link maps} both in work of Munson \cite{MunsonMfdCalcLink} and in our recent work \cite{CKKS}.  Thus $\mu_{123}$ should also appear in the (1,1,1) stage of the tower for the \emph{embedding} space $\LL_3^3$.  This would also match Conjecture 1.1 of \cite{BCSS}, which says that every type-$n$ invariant factors through the $(n+1)$-th stage of the homotopy tower.  (The triple linking number for string links is a finite-type invariant of type 2, and $1+1+1=3$.)  However, while stage (2,2,2) may not be the lowest stage at which $\mu_{123}$ appears, it does seem to be the best result that one can prove using any method based on configuration space integrals.  This is related to the reason why \cite{VolicConfig} yields type-$n$ invariants at stage $2n$, rather than $n+1$ as in Conjecture 1.1 of \cite{BCSS}: to construct a type-$n$ invariant via Bott--Taubes integrals, one needs $n$ \emph{pairs} of points, rather than just $n$ points.

\subsection{Organization of the paper}
In Section \ref{BTIs}, we briefly review finite-type invariants of knots and links.  We then review the configuration space integrals of Bott and Taubes, including the generalization to string links.  We also review our previous results concerning the homotopy-theoretic reformulation of these integrals.

In Section \ref{Tower}, we review the ``punctured knots" model for the stages of the Taylor tower, as well as Voli\'{c}'s result extending the Bott--Taubes integrals to the stages of the tower in this model.  We then show that this result generalizes to the setting of string links by proving Theorem \ref{FTLinkInvtsFactor}.  Thus finite-type invariants of string links factor through the Taylor tower, just as in the case of knots.  Finally, we show that the constructions of \cite{Rbo, HoMilnor3ple} extend to the Taylor tower, proving the main results, namely Proposition \ref{RboOnTaylor} and Theorem \ref{HoMilnorOnTaylor}.

In Section \ref{AlignedMaps}, we discuss an alternative method for obtaining the main results.  The alternative is to use a variant of the aligned maps model rather than the punctured knots model.  The advantage of the aligned maps model is that it has a monoid structure compatible with connect-sum, whereas we do not know of such a structure on the punctured knots model.  
With the appropriate aligned maps model, it is straightforward to carry out the main constructions in this paper (including the original result of Voli\'{c} in \cite{VolicConfig}).  Thus this model yields alternative proofs for all of these results about the tower.  We keep this last Section brief: the results we outline there are equivalent to those in the previous Section, and the potential benefits of combining them with the monoid structure remains to be explored.

\subsection{Acknowledgments}
This paper began as part of the author's Stanford University 2010 Ph.D. thesis under the direction of Ralph Cohen.  He is grateful to Ralph Cohen for his optimism, advice, and encouragement.  He thanks Tom Goodwillie for helping him find a crucial error in an earlier draft of this paper.  He thanks Ismar Voli\'{c} and Brian Munson for conversations related to the content of this paper.  He thanks referees for careful readings of this paper and useful comments.

\section{Finite-type invariants, Bott--Taubes integrals, and generalizations}
\label{BTIs}

\subsection{Finite-type invariants of knots and links}
\label{VassilievIntro}
The configuration space integrals of Bott and Taubes are closely related to finite-type invariants of knots and links, which can be defined in elementary terms.  We review the main ideas here, but refer the reader to, for example, Bar-Natan's paper \cite{BarNatan} for details.  When we say ``links'' in this Subsection, we include knots as a special case.  We also use ``links'' to mean either long links or closed links.

To give the definition, we first note that any $\R$-valued (or abelian-group-valued) invariant of oriented links can be extended to links with finitely many transversely self-intersections by inductively applying the ``Vassiliev skein relation'':
\[
V\left(\raisebox{-1.2pc}{\includegraphics[height=3.7pc]{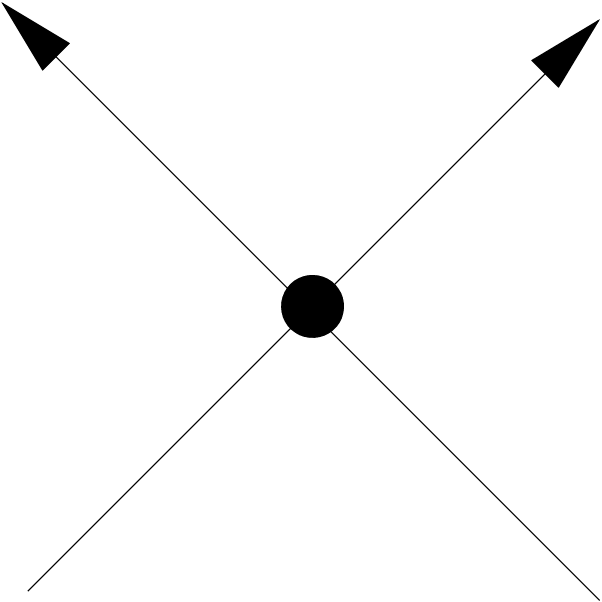}}\right) = V\left(\raisebox{-1.3pc}{\includegraphics[height=3.7pc]{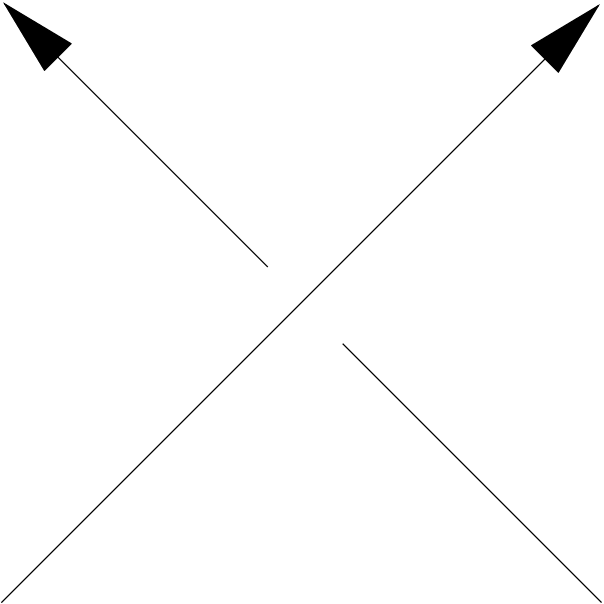}}\right) - V\left(\raisebox{-1.3pc}{\includegraphics[height=3.7pc]{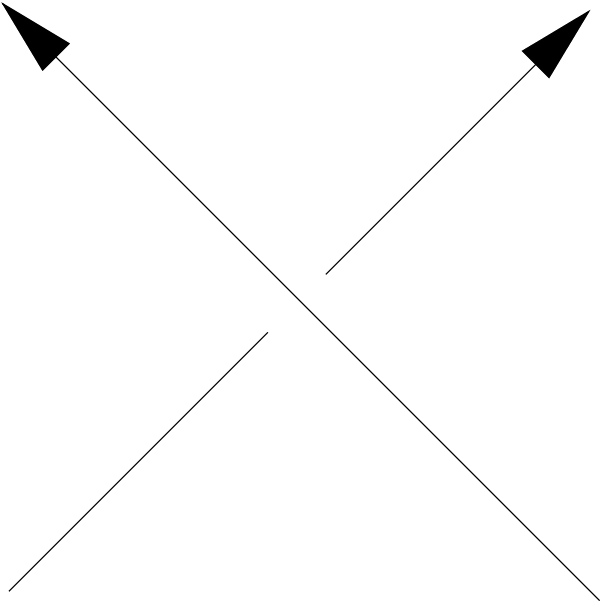}}\right)
\]
An invariant $v$ is then called \emph{finite-type} or \emph{Vassiliev} of type $n$ if $v$ vanishes on all knots with $>n$ self-intersections.  If we let $\V_n$ denote the $\R$-vector space (or $\Z$-module) of type-$n$ invariants, we have an increasing filtration $\V_0 \subset \V_1 \subset... \subset \V_n \subset ...$.  The conjecture of Vassiliev that any two knots can be distinguished by some finite-type invariant is still open.  

We next mention some of the simplest examples of finite-type invariants.  The space $\V_0$ is easily seen to be the space of constant functions.  The pairwise linking number (of either long links or closed links) is a type-1 invariant.  The triple linking number \emph{for long links}, which we write as $\mu_{123}$, is a type-2 invariant.  This invariant will feature prominently in Section \ref{HoMilnorOnTaylorSection}.

Our work in that Section relies on a connection between these invariants and certain combinatorial diagrams.
First, it is not difficult to construct a canonical, injective map from $\V_n / \V_{n-1}$ to a space $(\CD_n)^*$, which is the dual of a vector space of chord diagrams with $n$ chords.  Below is an element in $\CD_5$ in the setting of closed knots.  In the setting of long knots, the circle is replaced by an interval, and in the setting of closed (respectively long) $m$-component links, one has $m$ circles (respectively intervals). 
\[
\includegraphics[height=6pc]{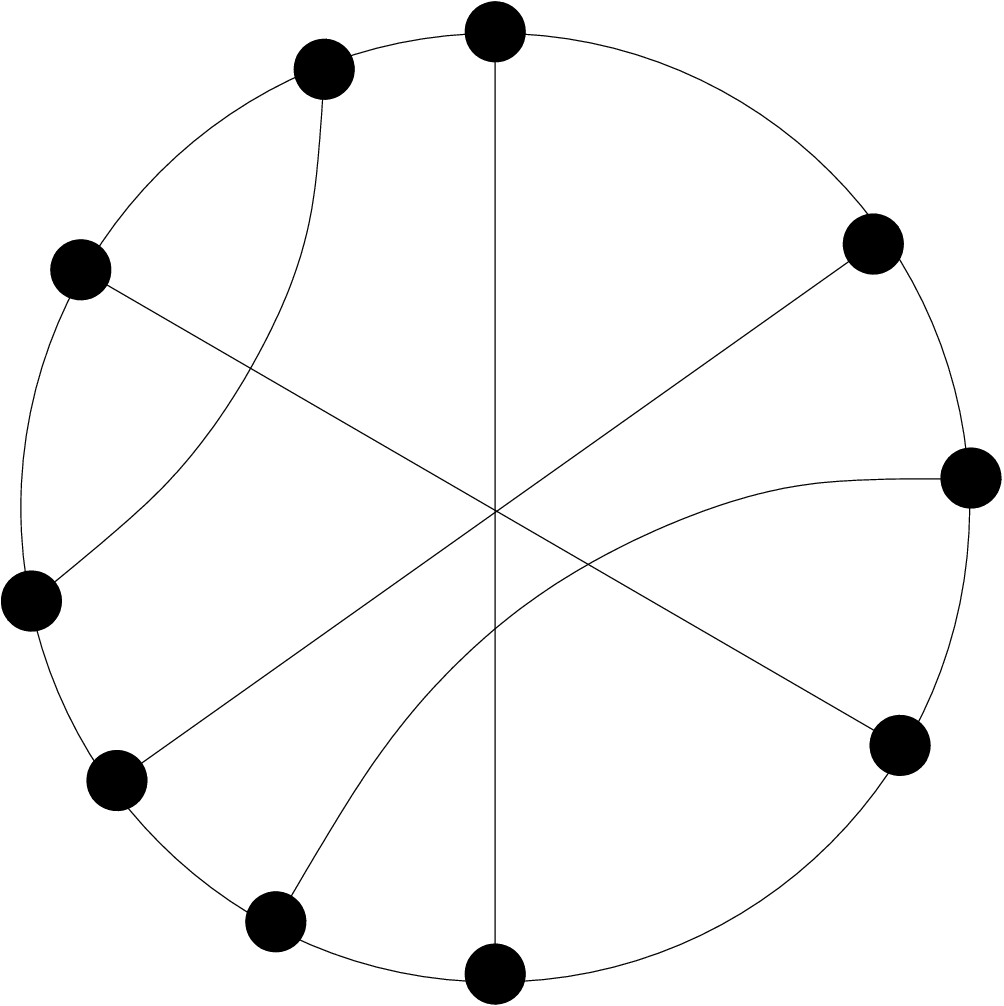}
\]
A much more difficult theorem is that this canonical map to $(\CD_n)^*$ is an isomorphism onto a certain subspace of $(\CD_n)^*$:
\[
\xymatrix{
W: \V_n / \V_{n-1}  \ar[r]^-{\cong} & (\CD_n / (1T,4T))^*}
\]
Above, $1T$ and $4T$ are the one-term and four-term relations imposed on chord diagrams in $\CD_n$; their definitions can be found in \cite{BarNatan}.  
The fact that this map $W$ is an isomorphism is deduced by constructing an inverse.  This inverse map can be constructed via either the Kontsevich integral or the Bott--Taubes integrals, the latter of which we discuss below.  We first review compactifications of configuration spaces, which are needed for Bott--Taubes integrals.

\subsection{Configuration spaces and their compactifications} 
\label{Compactifications}
The configuration space $C_n(M)$ is defined as the complement of the ``fat diagonal" in $M^n$: 
$$C_n(M):=\{(x_1,...,x_n) \in M^n | x_i \neq x_j \, \forall i \neq j\}.$$
It is not compact, even if $M$ is compact.  For a compact manifold $M$, there is a compactification of   Axelrod--Singer/Fulton--MacPherson (also sometimes called the ``canonical compactification"), which we denote $C_n[M]$.  This space keeps track of directions of collision and relative rates of approach of points in the configuration.  It is homotopy equivalent to the original configuration space $C_n(M)$.
It also has the important feature of being a smooth manifold with corners \cite{Axelrod-Singer, Fulton-MacPherson}.  This is crucial for the Bott--Taubes integrals, as well as for neat embeddings in our homotopy-theoretic reformulation.  
For $M=\R^d$
we define $C_n[\R^d]$ as a subspace of $C_{n+1}[S^d]$ with the last point at $\infty$.  
For $M=\R$, $C_n[\R]$ has $n!$ homeomorphic connected components; from now on, we will
consider just the component where the $n$ points are in order.

The compactification $C_n[M]$ can be defined via blowups.  Specifically, $C_n[M]$ is the closure of the image of 
\[
C_n(M) \incl \prod_{\underset{|S| \geq 2}{S \subset \{1,...,n\}}} \Bl(M^n, \Delta_S)
\]
where $\Bl(M^n, \Delta_S)$ denotes the blowup of $M^n$ along the diagonal where the points indexed by $S$ are equal.
This endows $C_n[M]$ with a stratification, where each stratum is indexed by a collection $\{S_1,...,S_k\}$ of (distinct) subsets $S_i \subset \{1,...,n\}$.  
The sets $\{S_1,...,S_k\}$ indexing strata are precisely those which are either nested or disjoint: $S_i \cap S_j \neq \emptyset$ implies either $S_i \subset S_j$ or $S_j \subset S_i$.  Intuitively, each subset $S_i$ indicates a collision of all the points whose indices are in $S_i$.  If $S_i \subset S_j$, then the points indexed by $S_i$ have collided faster than the remaining points indexed by $S_j$.  We note that a stratum indexed by $\{S_1,...,S_k\}$ has codimension $k$.
Finally, notice that taking $C_n[\R]$ to be just the component where the points are in order puts further restrictions on which $S_i$ can appear in a collection $\{S_1,...,S_k\}$ indexing a stratum of $C_n[\R]$.

  
There is another compactification of $C_n(M)$ called the simplicial compactification $C_n\la M \ra$.  This compactification is not smooth, but it is needed for the cosimplicial models for spaces of knots and links.  It keeps track of directions of collisions, but not relative rates of approach.  Thus $C_n\la M\ra$ is a quotient of $C_n[M]$.  (The reader may compare the definitions of the two compactifications in \cite[Definition 4.1]{SinhaTop}.)
The quotient map is a homotopy equivalence \cite[Corollary 5.9]{SinhaCptn}.  Define $C_n\la \R^d\ra$ as the subspace of $C_{n+1}\la S^d\ra$ where the last point is at $\infty$. 
As defined, $C_n\la \R \ra$ has $n!$ homeomorphic connected components; from now on, we will use $C_n\la \R \ra$ to denote just the component where the $n$ points are in order.  
Then $C_n\la \R \ra \cong \Delta^n$.  For us this is the only important point: there is a quotient map $C_n[\R] \to \Delta^n$ that just forgets relative rates of approach.

%

Finally, we will need one more slightly different compactification.  Let $C_n[I^d]$ and $C_n\la I^d \ra$ respectively denote the subspaces of $C_{n+2}[\R^d]$ and $C_{n+2}\la \R^d \ra$ where all the points are in the cube $I^d=[-1,1]^d$ and the first and last points are fixed at $(\pm 1,0,...,0)$.  
(These spaces are denoted $C_n[I^d,  d]$ and $C_n\la I^d, \d \ra$ in the papers of Sinha and \cite{BCSS, BCKS}.)
Then $C_n\la I \ra \cong C_n \la \R \ra \cong \Delta^n$.  
The space $C_n[I]$ is the $n$-dimensional associahedron $K_{n+2}$, which is a quotient of $C_n[\R]$, since, as defined above, the latter space records relative rates of approach to $\infty$ and not just relative rates of approach to $+\infty$ and $-\infty$ separately.  
It is well known that the strata of $K_{n+2}$ can be indexed by (partial) parenthetizations of $n+2$ letters.  We can see that this corresponds to the Fulton--Macpherson indexing of strata of $C_n[I]$ by sets of subsets.
In fact, call a subset $S$ of an ordered set $T$ \emph{consecutive} if for every $i,j,k\in T$ with $i<j<k$ and $i, k \in S$, we have that $j \in S$.  
Then each stratum of $C_n[\R]$ is indexed precisely by collections $\{S_1,...,S_k\}$ of consecutive proper subsets of $\{-\infty, 1,2,...,n, +\infty\}$ which are either nested or disjoint.
Such a collection is easily seen to correspond to a parenthetization.

\subsection{Bott--Taubes integrals}
\label{BTIsubsection}
In \cite{Bott-Taubes}, Bott and Taubes constructed knot invariants by considering a bundle 
\[
 \xymatrix{
F[q;t] \ar[r] & E[q;t]\ar[d] \\
 & \Emb(S^1,\R^3)}
\]
over the space of knots in $\R^3$.  The fiber $F[q;t]$ over a knot $K$ is a compactification of a configuration space of $q+t$ points in $\R^3$, $q$ of which lie in the image of $K$.  Another way of saying this is that the total space $E=E[q;t]$ is the pullback in the square below:
\begin{equation}
\label{BTsquare}
 \xymatrix{
E[q;t] \ar[r] \ar[d] \pullbackcorner & C_{q+t}[\R^3]\ar[d] \\
\Emb(S^1,\R^3) \x C_q[S^1] \ar[r] & C_q[\R^3]}
\end{equation}
Here the right-hand map is projection to the first $q$ points.  The bottom map is given by evaluating the embedding on the $q$ points in $S^1$.  Explicitly, this bottom map sends a point $(K,t_1,...,t_q)$ to $(K(t_1),...,K(t_q))$.

The cohomology of the configuration space $C_{q+t}[\R^3]$ is generated by certain ``spherical forms'' $\theta_{ij}$; each $\theta_{ij}$ is the pullback of the volume form on $S^2$ via the map $\phi_{ij}$ which records the unit vector between $x_i$ and $x_j$.
Bott and Taubes then considered the pullbacks of these $\theta_{ij}$ to $E[q;t]$, which by abuse of notation we will also call $\theta_{ij}$.  They integrated certain sums of products of these $\theta_{ij}$ along the fiber $F[q;t]$ of the bundle $E[q;t] \to \Emb(S^1, \R^3)$.  
Since $F[q;t]$ has nonempty boundary, Stokes' Theorem implies that 
\[
d \int_{F[q;t]} \alpha = \int_{\d F[q;t]} \alpha|_{\d F[q;t]}.
\]
Thus, to show that a fiberwise integral produces a closed form, it suffices to show that the integral along the boundary of the fiber vanishes.  Bott and Taubes showed that the integrals along certain types of boundary faces vanish, while for appropriate choices of $\alpha$, the contributions along the remaining faces cancel each other.  The latter type of face, along which the integrals do not vanish, is a \emph{principal face}, which is defined as any (codimension-one) stratum of $C_n[M]$ involving a collision of only two points,\footnote{Note that because of the blowups performed in constructing $C_n[M]$, any stratum indexed by a single set $\{S\}$ has codimension one, regardless of the cardinality of $S$.} i.e., any stratum indexed by $\{S\}$ with $S \subset \{1,...,n\}, |S|=2$.
The fact that principal face contributions do not vanish is important for recasting this construction in homotopy-theoretic terms.
Ultimately, Bott and Taubes produced a zero-dimensional closed form on $\Emb(S^1, \R^3)$, which represents a knot invariant.  

Their methods were then used by D.~Thurston to obtain all finite-type knot invariants \cite{Thurston, VolicBT} by constructing (for each $n$) an inverse to the map 
\begin{equation}
\label{Wmap}
W: \V_n/\V_{n-1} \to (\CD_n/(1T, 4T))^*
\end{equation}
mentioned in Section \ref{VassilievIntro}.
This requires enlarging the space $\CD_n$ of chord diagrams on $2n$ vertices to a space $\TD_n$ of trivalent diagrams on $2n$ vertices.  The picture below shows an element in $\TD_5$ in the setting of closed knots.
\[
\includegraphics[height=8pc]{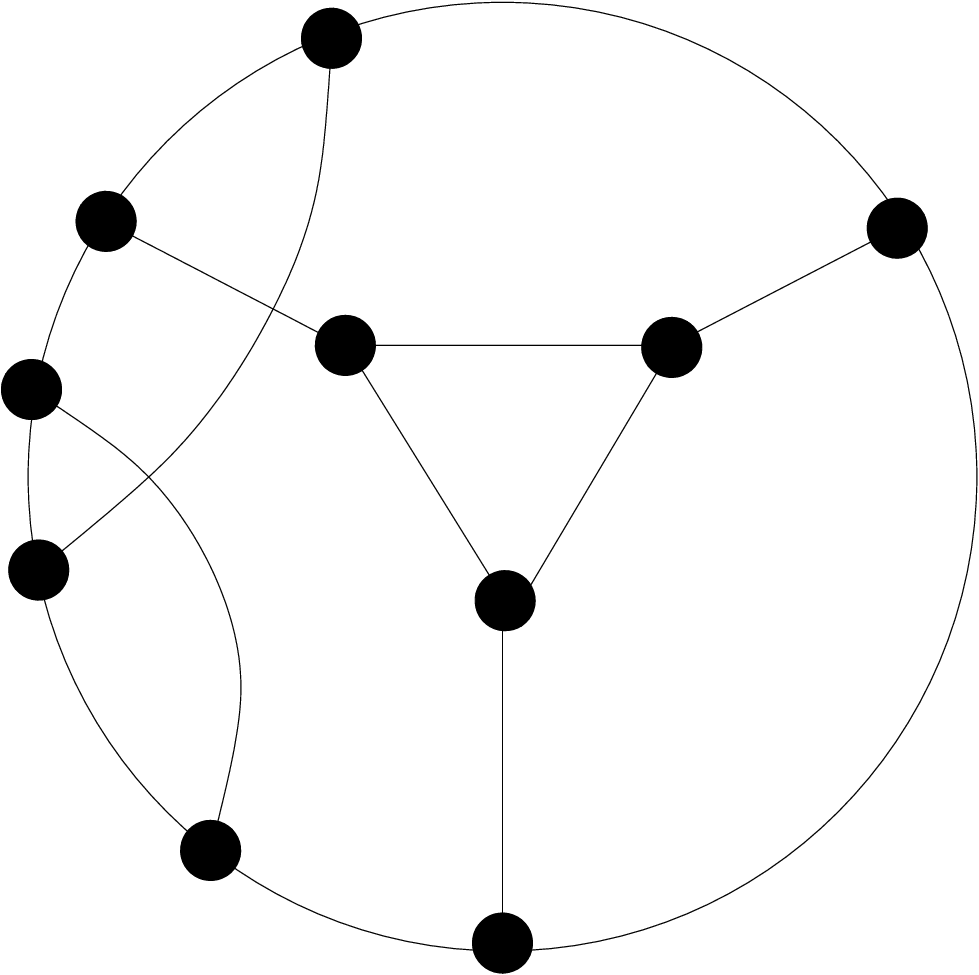}
\]
To a diagram $D \in \TD_n$ one associates the configuration space bundle $E_D=E[q;t]$ with fiber $F_D$, where $q$ is the number of vertices on strands and $t$ is the number of remaining vertices.  One also associates to $D$ a differential form $\theta_D = \bigwedge_{(i,j)\in D} \theta_{ij}$, where this wedge product is taken over all edges $(i,j)$ in $D$.  Thus a diagram $D$ produces a fiberwise integral $\int_{F_D} \theta_D$.  Counting the dimensions of $\theta_D$ and $F_D$ using the trivalence of $D$ shows that the result is a 0-dimensional form on $\Emb(S^1,\R^3)$.  
If one imposes a relation on $\TD_k$ called the STU relation, one finds that $\TD_k/(1T,STU) \cong \CD_k/(1T, 4T)$ \cite{BarNatan}.  
Thus an element $w \in \CD_k/(1T, 4T)$ corresponds to an element $w \in \TD_k/(1T,STU)$, and the inverse to the map $W$ in (\ref{Wmap}) is essentially given by\footnote{We have omitted the so-called \emph{anomaly term} in our formula for $BT$.  This (possibly nonzero) anomaly term is needed to construct finite-type knot invariants, but it vanishes in the case of finite-type link-homotopy invariants, such as $\mu_{123}$.} 
\begin{equation}
\label{InverseToW}
BT: w \mapsto BT(w)=\sum_{D \in \TD_n} w(D) \int_{F_D} \theta_D.
\end{equation}  
\begin{theorem}[\cite{Thurston, VolicBT}]
\label{BTIsAreUniversal}
Given any $\R$-valued type-$n$ knot invariant
\[
v: \pi_0(\Emb(S^1, \R^3)) \to \R
\]
there is an element $w \in \TD_k/(1T,STU)$ such that $v$ is given by the sum of configuration space integrals $BT(w)$:
\[
v=BT(w): \pi_0(\Emb(S^1, \R^3)) \to \R.
\]  
Explicitly, $BT(w)$ is given by (\ref{InverseToW}), and since $BT$ is inverse to the canonical map $W$, $w=W(v)$.
\end{theorem}

The square (\ref{BTsquare}) above can be generalized in several ways.  First, it is completely straightforward to replace closed knots by long knots: one just replaces $S^1$ by $\R$ in (\ref{BTsquare}).   Second, one can generalize from knots to links (long or closed).  The obvious analogue of Theorem \ref{BTIsAreUniversal} holds in all of these cases too.  
The generalization for long links, which is slightly less straightforward than that for closed links, is \cite[Theorem 5.6]{MunsonVolicHtpyLinks}.
The square (\ref{BTsquare}) also makes sense if one replaces $\R^3$ by $\R^d$, $d\geq 4$.  Cattaneo, Cotta-Ramusino, and Longoni constructed nontrivial real cohomology classes in $\Emb(S^1,\R^d)$ in this way \cite{Cattaneo}.  

The generalization to long links
requires one to keep track of relative rates of approach to infinity of configuration points in the domain.  In this most general setting of $\LL_m^d$, one considers the pullback
\begin{equation}
\label{LinkBTSquare}
\xymatrix{
E[q_1,...,q_m; t] \ar[r] \ar[d] \pullbackcorner & C_{q_1+...+q_m+t}[\R^d] \ar[d] \\
 \LL_m^d \x C_{q_1+...+q_m}[\coprod_{i=1}^m \R] \ar[r] & C_{q_1+...+q_m}[\R^d]  
}
\end{equation}
where the compactified configuration space lower-left corner is defined as follows:
\begin{definition}
\label{ClosureStringLinkImage}
Let $C_{q_1+...+q_m}[\coprod_{i=1}^m \R]$ be the closure of the image of the map
\begin{equation}
\label{StringLinkInducedMap}
C_{q_1}[\R] \x ... \x C_{q_m} [\R] \to C_{q_1+...+q_m}[\R^d]
\end{equation}
induced by any $m$-component string link.  
\end{definition}

This space was shown to be a manifold with corners in \cite[Lemma 4.4]{MunsonVolicHtpyLinks} and \cite[Section 2.3.2]{HoMilnor3ple}.

We now give another description of $C_{q_1+...+q_m}[\coprod_{i=1}^m \R]$ via blowups.  
This description will be useful in Section \ref{BTonTaylorLinks}.
First recall that the product $C_{q_1} [\R] \x ... \x C_{q_m} [\R]$ is a manifold with corners.  Recall also that a stratum in any factor $C_{q_i}[\R]$ is indexed by a collection of subsets of $\{1,...,q_i, \infty\}$, which indicates collisions of points.  
The strata in this product are of course just all possible products of strata in the factors.

\begin{proposition}
\label{BlowUpProp}
The space $C_{q_1+...+q_m}[\coprod_{i=1}^m \R]$ can be obtained from $C_{q_1} [\R] \x ... \x C_{q_m} [\R]$  by blowing up every stratum that involves collisions with $\infty$ on at least two of the $m$ strands.
\end{proposition}
\begin{proof}
First notice that if any two strata satisfy this property, then so does their intersection, so this blowup is well-defined, by first blowing up the lowest-dimensional strata.  
Next, in Definition \ref{ClosureStringLinkImage}, 
the map (\ref{StringLinkInducedMap}) is injective away from configurations in which some point is at $\infty$.  So away from $\infty$, $C_{q_1+...+q_m}[\coprod_{i=1}^m \R]$ and $C_{q_1} [\R] \x ... \x C_{q_m} [\R]$ agree via a stratum-preserving diffeomorphism.

Near $\infty$, $C_{q_1+...+q_m}[\coprod_{i=1}^m \R]$ acquires the stratification on $C_{q_1+...+q_m}[\R^d]$.  In particular, each codimension-$k$ stratum is given by a set of nested or disjoint sets 
\begin{align}
\label{StrataAtInfty1}
& \{S_1,...,S_k\} 
\end{align} 
\begin{align*} 
\mbox{ with each }  & S_i \subset \{1,...,q_1+...+q_m, \infty\}.
\end{align*} 
(The only difference between the two stratifications is that  points in a configuration in $C_{q_1+...+q_m}[\coprod_{i=1}^m \R]$ must lie in order on the string link strands, so not every stratum in $C_{q_1+...+q_m}[\R^d]$ occurs as a stratum in $C_{q_1+...+q_m}[\coprod_{i=1}^m \R]$.  Alternatively, we can think of such strata as empty in $C_{q_1+...+q_m}[\coprod_{i=1}^m \R]$.  Note also there are two directions of approach to infinity along a string link, so each stratum at infinity in $C_{q_1+...+q_m}[\coprod_1^m \R]$ is not connected, whereas each stratum at infinity in $C_{q_1+...+q_m}[\R^d]$ is connected.)  

On the other hand, a codimension-$k$ stratum at $\infty$ in $C_{q_1} [\R] \x ... \x C_{q_m} [\R]$ is given by
\begin{align}
\label{StrataAtInfty2}
(\{S^1_1,...,S^1_{k_1}\}, ..., \{S^m_1,...,S^m_{k_m}\}) \cong
\{S^1_1,...,S^1_{k_1}, ... , S^m_1,...,S^m_{k_m}\} 
\end{align}
\begin{align*}
\mbox{ with each } S^j_i \subset \{1+\sum_{\ell=1}^{j-1} q_\ell,...,\sum_{\ell=1}^j q_\ell, \infty\} \mbox{ and with } k_1+...+k_m=k.
\end{align*}
The difference between (\ref{StrataAtInfty1}) and (\ref{StrataAtInfty2}) is that in (\ref{StrataAtInfty1}), an $S_i$ may contain points from different strands.  Thus in the new stratification at $\infty$ (\ref{StrataAtInfty1}), we record relative rates of approach to $\infty$ of points on different strands (rather than just relative rates of approach of points on the same strand).  But the relative rates of approach of each such collection of points $\{x_1,...,x_p, \infty\}$ are recorded precisely by blowing up the stratum where the points $x_1,...,x_p, \infty$ have collided.
\end{proof}

\subsection{A related homotopy-theoretic construction} 
In previous work, we considered the Bott--Taubes bundle described above, and carried out ``integration along the fiber" homotopy-theoretically, which we briefly review now.   See our work \cite[Section 3]{Rbo} for details.  First, take a neat embedding
 of the total space $E=E[q;t]$ into a trivial bundle of a Euclidean space with corners:
 \[
 e_N: E \incl \K \x \R^{N-M} \x [0,\infty)^M.
 \]
Then collapse by the complement of a tubular neighborhood of $E$.  Quotienting by boundary subspaces gives a map
\[
\tau: \Sigma^N \K_+ \to E^{\nu_N}/(\d E^{\nu_N})
\]
from the $N$-fold suspension of the base space $\K$ (union a disjoint basepoint) to the Thom space  of the normal bundle $\nu_N$ of $e_N$, modulo its boundary.  In cohomology, this gives a map corresponding to integration along the fiber \cite[Corollary 3.7]{Rbo}.  

By letting $N$ in $e_N$ approach $\infty$, we then get a map from the suspension \emph{spectrum} of $\K$ to the Thom \emph{spectrum} of the normal bundle to the total space, which induces in cohomology a map similar to the Bott--Taubes integration along the fiber.   
Using the Thom isomorphism, this induces an ``integration along the fiber'' map in cohomology with arbitrary coefficients, producing classes in $H^*\K$.

\subsection{Recovering the triple linking number for string links}
\label{RecoveringTripleLinkingNumber}
In \cite{HoMilnor3ple}, we showed that a modification of the construction in \cite{Rbo} recovers $\mu_{123}$, the triple linking number for string links.  This modified construction is based on the fact that this invariant can be expressed as a sum of configuration space integrals \cite{MunsonVolicHtpyLinks}.  In \cite{HoMilnor3ple}, we identified the specific integrals in this sum, which correspond to four trivalent diagrams $L, M, R, T$, shown below.  A diagram gives rise to a configuration space integral, as explained in the case of knots in Section \ref{BTIsubsection}.
In the case of links, a diagram $D$ gives rise to the configuration space bundle $E_D =E[q_1,...,q_m;t] \to \LL_m^3$, where $q_i$ is the number of vertices on the $i$-th strand and $t$ is the number of vertices not on any strand, such as the one in the diagram $T$ below.  
\begin{figure}[h!]
\includegraphics[height=2.5pc]{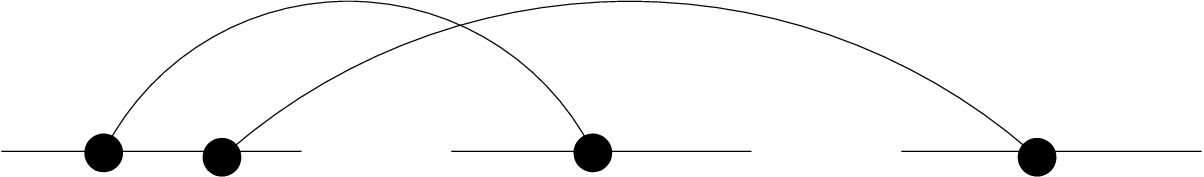} \qquad 
\includegraphics[height=3pc]{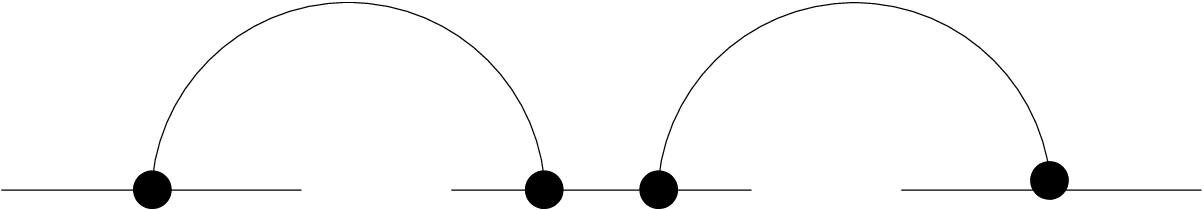} \\
\includegraphics[height=2.5pc]{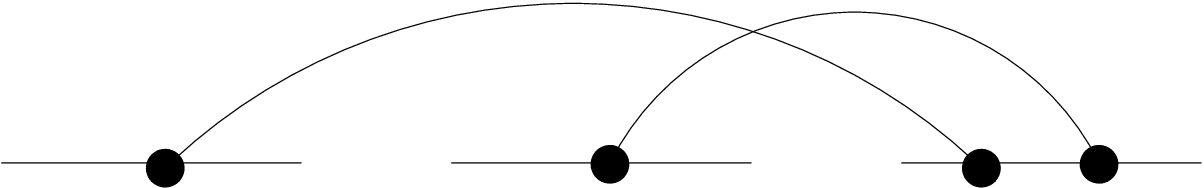} \qquad
\includegraphics[height=4.5pc]{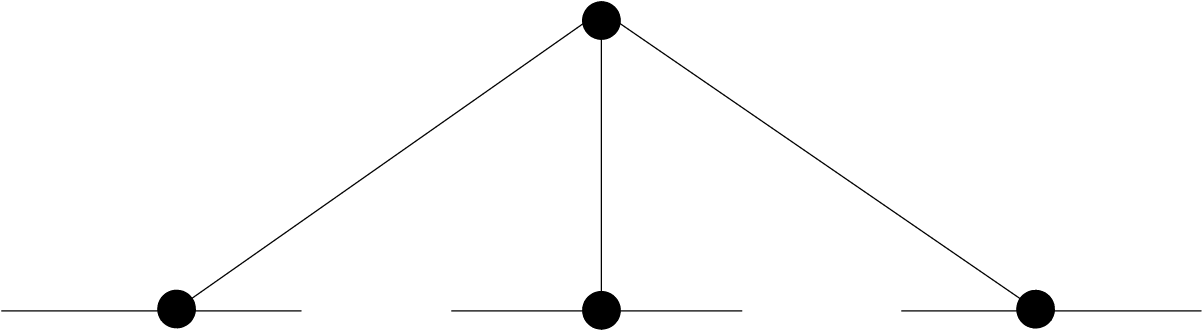}
\caption{The trivalent diagrams $L, M$ (top row), $R$, and $T$ (bottom row).}
\label{4UnitrivalentDiagrams}
\end{figure}
The four bundles $E_L \x S^2$, $E_M \x S^2$, $E_R \x S^2$, and $E_T$ thus all have 6-dimensional fibers.  We glued these bundles fiberwise along principal faces, according to the principal face cancellations in the integrals.

The resulting bundle is denoted $E_g \to \LL_3^3$.  The remaining boundary $\d E_g$ consists of faces at infinity.  
By gluing together neat embeddings of each bundle into $\LL_3^3 \x \R^M \x [0,\infty)^N$, we obtained a neat embedding of $E_g$ into a trivial bundle of glued Euclidean spaces with corners: 
\[
e: E_g \incl \LL_3 \x \left( \R^M \x [0,\infty)^N \right)/\sim.
\]  
We then took the Pontrjagin--Thom collapse and quotiented the remaining boundary to get a map in cohomology $\tau^*: H^*(E_g, \d E_g) \to H^{*-6}(\LL_3^3)$.  

To recover the triple linking number $\mu_{123}$, we used 
a map
\[
E_g \to S^2 \x S^2 \x S^2.
\]
On $E_T \subset E_g$, each of the three maps to $S^2$ corresponds to an edge in $T$, i.e., the map is the unit vector between the endpoints of the edge. 
On each remaining piece $E_D \x S^2 \subset E_g$, two of the three maps to $S^2$ are similar ``unit vector maps'' $E_D \to S^2$, while the third is given by projection to the $S^2$ factor.  The above map descends to a map of pairs
\[
(E_g, \d E_g) \to (S^2 \x S^2 \x S^2, \mathcal{D})
\]
where $\mathcal{D}$ is a subspace containing the image of $\d E_g$ (the faces at infinity).
We then considered a certain cohomology class $[\alpha] \in H^6(S^2 \x S^2 \x S^2, \mathcal{D})$ which maps to a generator of $H^6(S^2 \x S^2 \x S^2)$.  Let $[\beta]\in H^6(E_g, \d E_g)$ denote the pullback of $[\alpha]$ to $(E_g, \d E_g)$.  
One main result in \cite{HoMilnor3ple} is that the image of $[\beta]\in H^*(E_g, \d E_g)$ under $\tau^*$ is $\mu_{123} \in H^0(\LL_3^3)$.  
This result requires choosing the correct lift $[\alpha]$ of a generator in $H^6(S^2 \x S^2 \x S^2)$ to $H^6(S^2 \x S^2 \x S^2, \mathcal{D})$.  This choice is specified and justified in Sections 3.4 and 3.5 of our paper \cite{HoMilnor3ple}.

\section{Integration and its homotopy-theoretic analogue on the Taylor tower}
\label{Tower}

In \cite{VolicConfig}, Voli\'{c} showed that the Bott--Taubes integrals can be done on the stages of the Taylor tower for the knot space $\K=\Emb(\R, \R^d)$.  In this section we first observe that Voli\'{c}'s extension applies equally well to $\LL_m^d$.  We then show that our homotopy-theoretic constructions can also be carried out on these spaces.

\subsection{The ``punctured knots" model for the Taylor tower for the knot space}

We begin by briefly reviewing the Taylor tower for the knot space, which is also described in papers of Voli\'{c} \cite{VolicFT, VolicConfig} and Sinha \cite{SinhaTop, SinhaOperads}.  It consists of spaces $T_n\K$ with maps $T_n\K \to T_{n-1}\K$ and compatible maps $\K \to T_n\K$.
In a certain ``mapping space'' model of $T_n\K$, the map $\K \to T_n \K$ is essentially an evaluation map. 
In this paper, we focus instead on the ``punctured knots'' model (used by Voli\'c in \cite{VolicConfig}), but by slight abuse of terminology, we will still call this map the ``evaluation map.''
To describe the latter model of $T_n\K$, fix some closed disjoint subintervals $A_0,...,A_n$ of $I$.  These are the ``punctures."  Let $\P[n]$ be the category whose objects are subsets of $[n]=\{0,1,...,n\}$ and whose morphisms are inclusions.  Let $\P_0[n]$ be the full subcategory of $\P[n]$ consisting of nonempty subsets.

\begin{definition}
Let $EC_n$ be the functor from $\P[n]$ to spaces which on objects is given by  
\[
S \mapsto E_S:=\Emb(\R - \bigcup_{i\in S} A_i, \R^d).
\]  
On morphisms, an inclusion $S \subset T$ is sent to the restriction map $E_S \to E_T$.  
Thus $EC_n$ is an $(n+1)$-cubical diagram of embedding spaces.  Define $T_n\K$ as the homotopy limit of the restriction of $EC_n$ to $\P_0[n]$.  In other words, $T_n\K$ is the homotopy limit of the associated \emph{punctured cubical diagram}.
\qed
\end{definition}

More concretely, this homotopy limit $T_n \K$ can be identified as the subspace  
\begin{equation}
\label{SubspaceOfProduct}
T_n \K \subset \prod_{\emptyset \neq S \subset\{0,...,n\} } \Map(\Delta^{|S|-1}, E_S)
\end{equation}
such that all squares constructed from these maps and inclusions $S\incl T$ of subsets of $\{0,1,...,n\}$ commute.
Since $\K$ is the initial object in $EC_n$ (and hence the homotopy limit of that diagram), the map $\K \to T_n\K$ is just the map from the homotopy limit of the whole cube to the homotopy limit of the restriction to $\P_0[n]$.  (Also note that for $n\geq 2$, $\K$ is the limit of the diagram $EC_n|\P_0[n]$.)  In concrete terms, $\K$ sits inside $T_n \K$ as constant families of punctured knots.
The maps $T_n \K \to T_{n-1}\K$ are the restrictions of homotopy limits induced by an inclusion $[n-1] \incl [n]$.  (Any two choices of this inclusion induce homotopic maps.)

\subsection{The ``punctured links" model for the multi-tower for spaces of links} 
For the space $\LL_m^d$, the embedding calculus gives rise to a \emph{multi-tower} with stages indexed by $m$-tuples of natural numbers, rather than by natural numbers.  In this Subsection, we follow the work of Munson and Voli\'{c} \cite[Section 4]{MunsonVolicMultivar}.  We view $\Emb(\coprod_{i=1}^m \R, \R^d)$ as the value of the functor $\Emb(-, \R^d)$ from $\mathcal{O}(\coprod_{i=1}^m \R) = \prod_{i=1}^m \mathcal{O}(\R)$.  Thus we have a multivariable functor of $m$ variables.  

Define each stage in the tower as follows.  
Given a multi-index $\vec{n}=(n_1,...,n_m)$, fix $n_i+1$ closed disjoint subintervals $A^i_0,..., A^i_{n_i} \subset I$ in the $i$-th copy of $\R$.  Now let $EC_{\vec{n}}$ be the functor from $\P[n_1] \x ... \x \P[n_m]$ to spaces given by 
\[
(S_1,...,S_m) \mapsto E_{S_1,...,S_m} := \Emb \left(\coprod_{i=1}^m \left( \R - \bigcup_{j \in S_i} A^i_j \right), \R^d \right)
\]
Thus $EC_{\vec{n}}$ is a multi-cubical diagram.  Define $T_{\vec{n}} \LL_m^d$ as the homotopy limit of the restriction of $EC_{\vec{n}}$ to $\P_0[n_1] \x... \x \P_0[n_m]$.  
Thus $T_{\vec{n}} \LL_m^d$ is a subspace of
\[
\prod_{(S_1,...,S_m) \in \P_0[n_1] \x...\x \P_0[n_m] } \Map(\Delta^{|S_1|-1} \x ... \x \Delta^{|S_m| -1}, E_{S_1,...,S_m}).
\]

\subsection{Bott--Taubes integrals over the Taylor tower for the space of knots}
\label{BTIsOverTaylor}
As in (\ref{SubspaceOfProduct}), $T_n \K$ is a subspace of $\prod_S \Map(\Delta^{|S|-1}, E_S)$, a product of spaces of continuous maps.  One can put a smooth structure on the embedding spaces $E_S$ and replace $\Map(\Delta^{|S|-1}, E_S)$ by the space of \emph{smooth} maps\footnote{The space $\Delta^{|S|-1}$ can be given a smooth structure via charts from Euclidean spaces with corners.  A smooth map on $\Delta^{|S|-1}$ just corresponds to a map which extends to an open subset of Euclidean space.} 
$\Delta^{|S|-1} \to E_S$.
The resulting subspace of $T_n \K$ is weakly homotopy equivalent to $T_n\K$, roughly because (families of) continuous maps are homotopic to (families of) smooth maps.  Cf. \cite[Definition 2.3]{VolicConfig}.
From now on, we will take $\Map(\Delta^{|S|-1}, E_S)$ and $T_n \K$ to mean spaces of smooth maps.  This allows us to put the structure of a smooth infinite-dimensional manifold $T_n \K$, much like one can for an embedding space such as $E_S$.  
We similarly take $T_{\vec{n}} \LL^d_m$ to be a space of smooth maps, which thus admits a smooth structure.  This is needed for constructing differential forms in these spaces via fiberwise integration.

Now one would like to define a bundle over $T_n\K$ by replacing the space of knots $\K$ by $T_n\K$ in  the original Bott--Taubes pullback square (\ref{BTsquare}).  Recall 
that the quotient $C_n\la I \ra$ of $C_n[I]$ is diffeomorphic as a manifold with corners to $\Delta^n$.  Thus we can include $C_n[I] \incl C_n[I] \x \Delta^n$ as the graph of the quotient map.  Using the inclusion (\ref{SubspaceOfProduct}), we then \emph{almost} have the composition 
\[
\xymatrix{
T_n\K \x C_n[I] 
\ar@{^(->}[r] & 
\left(\prod_{S\subset \{0,...,n\}} \Map(\D^{|S|-1},E_S) \x \Delta^n\right) \x C_n [I] \ar[r] & 
C_n[\R^d]
} 
\]
The problem is that each $E_S$ is a space of knots with $|S|$ punctures, and a point in a configuration in $C_n[I]$ may lie in one of the punctures.
Voli\'{c} resolves this difficulty in \cite{VolicConfig} by following the canonical quotient map $C_n[I] \to C_n\la I \ra \cong \Delta^n$ by a smooth map $\gamma : C_n\la I\ra \to \Delta^n$, which depends on the punctures $A_0, ..., A_n$.  Roughly, this map $\gamma$ sends a configuration with points lying in certain punctures to a certain boundary point $\vec{t}$; this point $\vec{t}$ is in turn sent by an element of holim$(EC_n)$ to a punctured knot which has all of those punctures filled in.

More precisely, write the target of $\gamma$ in barycentric coordinates: $\Delta^n = \{(t_0,...,t_n)| t_0+...+t_n=1\}$.  Each codimension-one face of $\Delta^n$ can be written as $\d_i\Delta^n:=\{(t_0,...,t_n)| t_i=0\}$ for $i=0,..,n$.
Fix a configuration $c$ in $C_n[I]$ with points located at $(x_1,...,x_n)$.  Let $T:=\{i \in \{0,...,n\} | \mbox{ some } x_j \in A_i\}$ and let $S=[n]\setminus T$.
Then $\gamma$ has the property that 
\begin{equation}
\label{GammaProperty1}
\gamma(c) \in \bigcap_{i \in T} \d_i \Delta^n. 
\end{equation}
Now an element of $T_n \K$ is a collection of maps, one of which is a map from $\Delta^{|S|-1} \cong \Delta^{n-|T|} \cong \bigcap_{i \in T} \d_i \Delta^n$ to $E_S = \Emb(\R - \sqcup_{i \notin T} A_i, \R^d)$.
Thus property (\ref{GammaProperty1}) ensures that for every configuration in $C_n[I]$ with points located at $(x_1,...,x_n)$, we get a point in some $\Delta^{|S|-1}$ which when plugged into an element of $\Map(\Delta^{|S|-1}, E_S)$ 
gives a punctured knot $e_S\in E_S$
such that no $x_j$ lies in any of its punctures.

Another key feature of $\gamma$ is that its restrictions to any two faces are the same map from an $(n-1)$-simplex (and the restrictions of that map agree on all the $(n-2)$-dimensional faces, etc.).
Thus for every $q\leq n$, we have the restriction $\gamma_q: C_q\la I \ra \to \Delta^n$.

The map $\gamma$ is constructed inductively on $n$ and via partitions of unity on the simplex, which ensure smoothness.  We refer the reader to Proposition 4.1 and its proof in \cite{VolicConfig} for all the details.

Voli\'{c} then defines the space $\Gamma[q]$ as the graph of the composite $C_q[I] \to C_q \la I \ra \to \Delta^n $.
The properties of $\gamma$ allow the bottom arrow in the square below to be defined for any $q \leq n$, and $\Gamma[q;t]$ is then defined as the pullback in this square.
\begin{equation}
\label{GammaBTSquare}
\xymatrix{ \Gamma[q;t] \pullbackcorner \ar[r] \ar[d] & C_{q+t}[\R^d] \ar[d] \\
\Gamma[q] \x T_n\K \ar[r] & C_q[\R^d]}
\end{equation}
Projecting to $T_n\K$ gives a bundle $\Gamma[q;t] \to T_n\K$ for any $n\geq q$, whose fiber is a smooth manifold with corners.  This allows one to integrate along the fiber just as in the original Bott--Taubes construction and produce invariants $T_n \K \to \R$ by configuration space integrals.  On the image of $\K$, this map agrees with the usual Bott--Taubes integrals.  This leads to part (1) of Theorem \ref{FTKnotInvtsFactorThruTower} below.
 
Part (2) of Theorem \ref{FTKnotInvtsFactorThruTower} concerns the stage of the tower at which one sees type-$n$ invariants and is deduced as follows.
As in Theorem \ref{BTIsAreUniversal}, a universal type-$n$ invariant requires configurations of $2n$ points, i.e., $q+t=2n$.  Moreover, one needs to consider configurations with all $q+t$ points on the knot, i.e., $t=0$ and $q=2n$ (which corresponds to chord diagrams).  Thus the lowest stage of the tower where one sees type-$n$ invariants is stage $2n$.

\begin{theorem}[Theorem 1.1 of \cite{VolicConfig}]
\label{FTKnotInvtsFactorThruTower}
(1)  Finite-type invariants of long knots over $\R$ factor through the Taylor tower $\{T_n \K\}_n$.  
(2) More specifically, any type-$n$ invariant factors through $T_{2n}\K$.  
In other words, given a type-$n$ invariant $v: \pi_0 \K \to \R$, there is a configuration space integral which gives the dashed arrow shown below and makes the diagram commute:
\[
\xymatrix{
 & \pi_0(T_{2n} \K) \ar@{-->}[dd] \\
\pi_0 \K \ar[ur]^-{\pi_0(ev_{2n}) \,\,\,\,\,\,\,\,} \ar[dr]_v & \\
 & \R
}
\]
\qed
\end{theorem}

\begin{remark}[Configurations on knots in a box vs. configurations on long knots]
\label{KnotsInABox}
Although Voli\'{c} seems to work entirely with configurations in $I$ in \cite{VolicConfig}, it seems preferable, or perhaps even necessary, to use configurations in $\R$, since one needs the vanishing of boundary integrals along faces at infinity in the Bott--Taubes construction.  Thus we replace $C_n[I]$ by $C_n[\R]$.  We again start with the canonical quotient $C_n[\R] \to C_n\la \R\ra$ by forgetting relative rates of approach.  We then follow this by the map $C_n \la \R \ra \to C_n \la I \ra$ given by sending points outside of $I$ to the nearest endpoint of $I$.  Finally, applying $\gamma$ gives a map $\widetilde{\gamma}: C_n [\R] \to \Delta^n$, defined as the composition
\[
\xymatrix{
\widetilde{\gamma}: C_n[\R] \ar[r] &  C_n \la \R\ra \ar[r] & C_n \la I \ra \ar[r]^-\gamma & \Delta^n.
}
\]
This map can easily be made smooth by smoothing out $\gamma$ near the boundary of $C_n \la I \ra$.  Again, there are restrictions $\widetilde{\gamma}_q: C_q [\R] \to \Delta^n$.  Thus from now on, we will take $\Gamma[q]$ to be the graph of $\widetilde{\gamma}$ and $\Gamma[q;t]$ to be the pullback in square (\ref{GammaBTSquare}) using this replacement of $\Gamma[q]$. 
\qed
\end{remark}

\subsection{Generalizing to link spaces}
\label{BTonTaylorLinks}
In this section, we will show that Voli\'{c}'s Theorem \ref{FTKnotInvtsFactorThruTower} generalizes from long knots to long links.
\begin{theorem}
\label{FTLinkInvtsFactor}
Finite-type invariants of string links over $\R$ factor through the Taylor tower $T_{\vec{\bullet}} \LL_m^3$ for any $m\geq 1$.  In general, an invariant of type $n$ factors through all stages at least as low as $(2n,2n,...,2n)$.  In other words given a type-$n$ invariant $v: \pi_0 \LL_m^3 \to \R$, there is a configuration space integral which gives the map shown by the dotted arrow below and makes the diagram commute:
\[
\xymatrix{
 & \pi_0(T_{(2n,2n,...,2n)}) \LL_m^3 \ar@{-->}[dd] \\
\pi_0 \LL_m^3 \ar[ur]^-{\pi_0(ev_{(2n,2n,...,2n)}) \,\,\,\,\,\,\,\,} \ar[dr]_v & \\
 & \R
}
\]
\end{theorem}

\begin{proof}
Fix $\vec{n}=(n_1,...,n_m)$.
We want to evaluate an element of $T_{\vec{n}} \LL_m^d$ on a configuration in $C_{q_1+...+q_m}[\coprod_{i=1}^m \R]$ where $q_i \leq n_i$ for all $i$.
For each $i$, find a map $\widetilde{\gamma}^{n_i} : C_{n_i} [\R] \to \Delta^{n_i}$ like $\widetilde{\gamma}$ in Remark \ref{KnotsInABox} above, and for any $q_i \leq n_i$, let $\widetilde{\gamma}^{n_i}_{q_i}$  denote the restriction of  $\widetilde{\gamma}^{n_i}$ to $C_{q_i}[\R]$.  By the construction of $\gamma$, this is the same map for any inclusion $\{1,...,q_i\} \incl \{1,...,n_i\}$.

By Proposition \ref{BlowUpProp}, $C_{q_1+...+q_m}\left[\coprod_{i=1}^m \R\right]$ is a blowup of $C_{q_1}[\R] \x ... \x C_{q_m}[\R]$, so there is a blow-down map 
\[
C_{q_1+...+q_m}\left[\coprod_{i=1}^m \R\right] \to C_{q_1}[\R] \x ... \x C_{q_m}[\R]
\]
given by forgetting relative rates of approach to $\infty$ of points on different strands.
Following this map by $\widetilde{\gamma}^{n_1}_{q_1} \x ... \x \widetilde{\gamma}^{n_m}_{q_m}$ gives a map
\begin{equation}
\label{GammaTilde}
\widetilde{\gamma}: C_{q_1+...+q_m}\left[\coprod_{i=1}^m \R\right] \to \Delta^{n_1} \x ... \x \Delta^{n_m}.
\end{equation}
If we define $\Gamma[q_1,...,q_m]$ as the graph of this map above, then the bottom horizontal arrow in the square below is well defined, and we can define $\Gamma[q_1,...,q_m; t]$ as the pullback in this square.  Note that $\Gamma[q_1,...,q_m] \cong C_{q_1+...+q_m}[\coprod_{i=1}^m \R] $.
\begin{equation}
\label{TaylorLinkBTSquare}
\xymatrix{ \Gamma[q_1,...,q_m;t] \ar[r] \ar[d] & C_{q_1+ ... +q_m+ t}[\R^d] \ar[d] \\
\Gamma[q_1,...,q_m] \x T_{\vec{n}}\LL_m^d \ar[r] & C_{q_1+ ... +q_m}[\R^d]}
\end{equation}
Following the left-hand vertical map by the projection to $T_{\vec{n}}\LL_m^d$ gives a map 
\[
\xymatrix{
\Gamma[q_1,...,q_m;t] \ar[d] \\ T_{\vec{n}}\LL_m^d
}
\]
over the $\vec{n}$-th Taylor tower stage for any $(q_1,...,q_m)$ with $q_i \leq n_i$ for all $i$.  

This map is a composition of two fiber bundles, where all the spaces involved are (possibly infinite-dimensional) manifolds with corners.  Hence it is a fiber bundle.  (This can be shown using Theorem 9.8 in Chapter 4 of Husemoller's book \cite{Husemoller}, since all the spaces involved are paracompact and the base of the composite bundle is locally contractible.)  

Since we want to integrate over its fiber, we now check that this fiber is a manifold with corners.  
By viewing $\Gamma[q_1,...,q_m;t] \to T_{\vec{n}}\LL_m^d$ as a composition of two bundles, we deduce that its fiber is the total space of a bundle with fiber and base given by the fibers of the two bundles.
That is, the fiber of $\Gamma[q_1,...,q_m;t] \to T_{\vec{n}}\LL_m^d$ is the total space $E$ of a bundle $F \to E \to \Gamma[q_1,...,q_m]$.
Here $F$ is the fiber of the left-hand vertical map in (\ref{TaylorLinkBTSquare}), which is the fiber of the right-hand vertical map in that square, which in turn is a manifold with corners.  The space $\Gamma[q_1,...,q_m] \cong C_{q_1+...+q_m}[\coprod_{i=1}^m \R]$ is also a manifold with corners.  
Hence $E$ is a manifold with corners, and we conclude that the fiber of this ``Bott--Taubes bundle over the Taylor tower'' is indeed a manifold with corners.  (In fact, as we will see in the next Subsection, it is isomorphic to the fiber of the usual Bott--Taubes construction over the link space.)

Finally, it is not difficult to see that the fiberwise integration over this bundle agrees with the usual Bott--Taubes integrals for long links on the image of $\LL^d_m \to T_{\vec{n}} \LL^d_m$.  This, combined with the long link analogue \cite[Theorem 5.6]{MunsonVolicHtpyLinks} of Theorem \ref{BTIsAreUniversal}, completes the proof.
\end{proof}

\begin{remark}
For a given finite-type invariant, we may be able to deduce the factoring at a lower stage than $(2n,2n,...,2n)$.  For example, by Theorem \ref{FTKnotInvtsFactorThruTower}, a type-$n$ knot invariant of the first strand factors through the stage $(2n,0,0,...,0)$.  But since an arbitrary type-$n$ invariant may require an integral over a configuration space with $2n$ points on any strand, this seems to be the best general bound that we can put on the multi-index $(n_1,...,n_m)$. \qed
\end{remark}

\subsection{Pontrjagin--Thom construction on the bundles over the tower}
\label{HoBTonTaylor}

We now extend the results of the previous Subsections to our homotopy-theoretic constructions from \cite{Rbo, HoMilnor3ple}.  To do this, we start by more carefully examining the bundles over the stages of the tower.


We start with (\ref{TaylorLinkBTSquare}) and replace $T_{\vec{n}} \LL_m^d$ by the image of the ``evaluation map'' $\LL_m^d \to T_{\vec{n}} \LL_m^d$:
\begin{equation}
\label{ImageOfEvBTSquare}
\xymatrix{ \Gamma[q_1,...,q_m;t]|_{\mathrm{im}(ev_{\vec{n}})} \ar[r] \ar[d] & C_{q_1+ ... +q_m+ t}[\R^d] \ar[d] \\
\Gamma[q_1,...,q_m] \x \mathrm{im}(ev_{\vec{n}}) \ar[r] & C_{q_1+ ... +q_m}[\R^d]}
\end{equation}
The evaluation map is injective, since it includes a link as a constant family of punctured links.  Thus its image can be identified with $\LL_m^d$, and the spaces and maps agree (up to diffeomorphism) with the ones in the original Bott--Taubes pullback square for long links (\ref{LinkBTSquare}).  
So if we consider the following commutative diagram of spaces, where the horizontal arrows are induced by the evaluation map, we see that those horizontal arrows are diffeomorphisms.
\[
\xymatrix{
E[q_1,...,q_m;t] \ar[r]^-\cong \ar[d] & \Gamma[q_1,...,q_m;t]|_{\mathrm{im}(ev_{\vec{n}})} \ar[d] \\
C_{q_1+...+q_m}\left[ \coprod_{i=1}^m \R \right] \x \LL_m^d \ar[r]^-\cong \ar[d] &
\Gamma[q_1,...,q_m] \x \mathrm{im}(ev_{\vec{n}}) \ar[d] \\
\LL_m^d \ar[r]^-\cong  & \mathrm{im}(ev_{\vec{n}})
}
\]
Hence over $\mathrm{im}(ev_{\vec{n}})$, this newly defined Bott--Taubes bundle is isomorphic to the original Bott--Taubes bundle.  Since we deduced in the previous Subsection that $\Gamma[q_1,...,q_m;t]$ is a locally trivial fiber bundle, we see that every fiber is isomorphic to the fiber in the original Bott--Taubes bundle.

We now want to embed these bundles into trivial bundles whose fibers are Euclidean spaces with corners.  We want these embeddings to be \emph{neat}, i.e., to preserve the corner structure.  
In \cite{Rbo} and \cite{HoMilnor3ple}, we made this precise via work of Laures \cite{Laures}.  Specifically, one considers \emph{manifolds with faces}, which, when equipped with $N$ codimension-one strata satisfying certain conditions, are called $\la N \ra$-\emph{manifolds}.  These are precisely the manifolds with corners which for some $M$ can be embedded neatly in $\R^M \x [0,\infty)^N$.  The main points for us are that compactified configuration spaces are $\la N \ra$-manifolds for some $N$ and that a neat embedding has a well defined tubular neighborhood (diffeomorphic to the normal bundle).
The reader may consult Section 2.1 (and especially Proposition 2.1.7) of \cite{Laures} or Section 3.1 of \cite{Rbo} for precise definitions and details.

Now if we give the total space $\Gamma[q_1,...,q_m;t]$ a corner structure from its fibers, we get the 
same corner structure as on $E[q_1,...,q_m;t]$.  
Carrying out neat embeddings $\Gamma[q_1,...,q_m;t] \incl T_{\vec{n}} \LL_m^d \x \R^M \x [0,\infty)^N$ and Thom collapse maps yields the following result, which says that our basic construction from \cite{Rbo} extends to the stages of the Taylor tower.  Its proof is similar to that of Theorem \ref{HoMilnorOnTaylor} in the next Subsection (which treats the construction with glued bundles in \cite{HoMilnor3ple}).  Since the construction in \cite{HoMilnor3ple} essentially supersedes the one in \cite{Rbo}, Theorem \ref{HoMilnorOnTaylor} is stronger than Proposition \ref{RboOnTaylor}.  Thus we omit the proof of Proposition \ref{RboOnTaylor} below and leave the interested reader to adapt the proof of Theorem \ref{HoMilnorOnTaylor}.

\begin{proposition}
\label{RboOnTaylor}
There is a commutative square of spectra as below on the left, where 
\begin{itemize}
\item
the superscript $\nu$ denotes the Thom space of the normal bundle of a certain embedding, 
\item
the vertical maps are ``evaluation maps," and 
\item
the horizontal maps are Thom collapse maps.
\end{itemize}
Applying cohomology, together with suspension isomorphisms and Thom isomorphisms, induces a commutative square in cohomology with arbitrary coefficients, as below on the right.
\[
\xymatrix{
\Sigma^\infty \K_+ \ar[r] \ar[d] & E[q;t]^\nu/ \d E[q;t]^\nu \ar[d] & & H^{*-\dim F[q;t]}(\K) & H^*(E[q;t], \d E[q;t]) \ar[l] \\
\Sigma^\infty (T_n\K)_+ \ar[r] & \Gamma[q;t]^\nu/ \d \Gamma[q;t]^\nu & & H^{*-\dim F[q;t]}(T_n\K) \ar[u] & H^* (\Gamma[q;t], \d \Gamma[q;t]) \ar[l] \ar[u]
}
\]
\end{proposition}

\subsection{Extending to the ``gluing" refinement and the triple linking number}
\label{HoMilnorOnTaylorSection}

We now consider the refinement in \cite{HoMilnor3ple} which recovers the triple linking number.  So fix $m=3, d=3$.  
Recall the trivalent diagrams $L,M,R,T$ shown in Figure \ref{4UnitrivalentDiagrams}, which are needed to construct $\mu_{123}$.  
In \cite{HoMilnor3ple}, we built a glued space $E_g$ out of the four spaces $E_L:=E[2,1,1;0], E_M:=E[1,2,1;0], E_R:=[1,1,2;0], E_T:=E[1,1,1;1]$.  More precisely, $E_g = (E_L \x S^2) \sqcup (E_M \x S^2) \sqcup (E_R \x S^2) \sqcup E_T / \sim$, where each $E_D \x S^2$ is glued to $E_T$ along a principal face.  

Using similar notation, we consider spaces $\Gamma_D$ for $D=L,M,R,T$ which fiber over $T_{\vec{n}}\LL_3^3$.  
Considering the indices $(q_1,q_2,q_3)$ for these four diagrams $D$, the ``smallest" $\vec{n}$ for which every $\Gamma_D$ defines a Bott--Taubes bundle over $T_{\vec{n}}\LL_3^3$ is $(2,2,2)$.  So for the rest of this Section, we consider only multi-indices $\vec{n}$ with $n_i \geq 2$ for all $i=1,2,3$.
We can then glue the spaces $\Gamma_L \x S^2 ,\Gamma_M \x S^2, \Gamma_R \x S^2$, and $\Gamma_T$ to create a space $\Gamma_g$ which fibers over $T_{\vec{n}}\LL_3^3$.  This is straightforward because the gluings are done fiberwise over a fixed base; the $E$-bundles and $\Gamma$-bundles have different bases, but isomorphic fibers.  

\begin{theorem}
\label{HoMilnorOnTaylor}
For any $\vec{n}=(n_1,n_2,n_3)$ with each $n_i\geq 2$, there is a commutative diagram in cohomology with arbitrary coefficients:
\[
\xymatrix{
H^{*-6}(\LL_3^3)  & H^*(E_g,\d E_g) \ar[l]  & H^*(S^2 \x S^2 \x S^2, \mathcal{D}) \ar[l] \ar[dl]\\
H^{*-6}(T_{\vec{n}}\LL_3^3) \ar[u] & H^* (\Gamma_g, \d \Gamma_g) \ar[l] \ar[u]  &
}
\]
Recall the class $[\beta] \in H^6(E_g, \d E_g)$ (see Section \ref{RecoveringTripleLinkingNumber})
whose image in $H^0(\LL_3^3)$ is $\mu_{123}$. 
This class $[\beta]$ can be lifted to $[\widetilde{\beta}] \in H^6(\Gamma_g, \d \Gamma_g)$.  Thus the homotopy-theoretic Bott--Taubes integrals provide a way of factoring $\mu_{123}$ through the Taylor tower for any stage as low as $T_{(2,2,2)}\LL_3^3$.
\end{theorem}

\begin{proof}
We want to construct a neat embedding and Pontrjagin--Thom collapse map for $\Gamma_g$ over the Taylor stage in such a way that these are compatible with the neat embedding and collapse map for the bundle over the link space itself.

In \cite[Lemma 7.1]{HoMilnor3ple}, we showed that there is a neat embedding 
\[
E_g \incl \LL_3^3 \x \left(\coprod_1^4 (\R^M \x [0,\infty)^N)/\sim\right)
\]
 for some $M,N$.  The main idea there is to start with embeddings of the various $E_D$ into $\LL_3^3 \x C_4[\R^3]$; neatly embed copies of that space into copies of $\LL_3^3 \x \R^M \x [0,\infty)^N$; glue those copies together; and then show that for $D=L,M,R$,  for sufficiently high-dimensional $\R^M$, the embedding of $S^2 \x$(boundary face of $E_D$) extends to an embedding of all of $S^2 \x E_D$.  We now extend this argument to $\Gamma_g$.  Notice first that any $\Gamma[q_1,...,q_m;t]$ embeds neatly into $T_{\vec{n}} \LL_m^d \x C_{q_1+...+q_m+t}[\R^d]$.  So for any $D\in \{L,M,R,T\}$, we have an embedding $\Gamma_D \incl
T_{\vec{n}}\LL_3^3 \x C_4[\R^3] $ such that the left-hand square below commutes.  The right-hand square, and its commutativity, come from neatly embedding $C_4[\R^3]$ into $\R^M \x [0, \infty)^N$.
\begin{equation}
\label{EmbedIntoConfig}
 \xymatrix{
E_D \ar@{^(->}[r] \ar[d] & 
\LL_3^3 \x C_4[\R^3]   \ar[d]  \ar@{^(->}[r] &
\LL_3^3 \x \R^M \x [0, \infty)^N \ar[d] \\
\Gamma_D \ar@{^(->}[r]  & 
T_{\vec{n}}\LL_3^3 \x C_4[\R^3]  \ar@{^(->}[r] &
T_{\vec{n}}\LL_3^3 \x \R^M \x [0, \infty)^N 
}
\end{equation}
Then for $D=L,M,R$, we claim that (just as for the $E_D$) we can embed $\Gamma_D \x S^2$ so that on a collar of the principal face boundary of $\Gamma_D$, this embedding can be glued to the embedding of $\Gamma_T$ along the appropriate principal face of $\Gamma_T$.  Away from the principal face boundary of $\Gamma_D$, the $S^2$ factor is embedded in a standard way.  We omit the details of such an argument, since they are the same as for $E_g$ in \cite[Lemma 7.1]{HoMilnor3ple}.  As a result, we get a neat embedding $\Gamma_g \incl T_{\vec{n}} \LL_3^3 \x \left(\coprod_1^4 (\R^M \x [0,\infty)^N)/\sim\right)$.  Since the incorporation of the $S^2$ factors away from the boundary of $\Gamma_T$ (respectively $E_T$) makes no use of the punctured links (respectively links), the following square commutes:
\begin{equation}
\label{EmbedIntoEuclid}
 \xymatrix{
E_g \ar@{^(->}[r] \ar[d] & 
\LL_3^3 \x \left(\coprod_1^4 (\R^M \x [0,\infty)^N)/\sim\right) \ar[d] \\
\Gamma_g \ar@{^(->}[r]  & 
T_{\vec{n}} \LL_3^3 \x \left(\coprod_1^4 (\R^M \x [0,\infty)^N)/\sim\right)
}
\end{equation}

We now want to collapse complements of tubular neighborhoods of the horizontal embeddings above in a compatible way.  For $\eps>0$, we will first define neighborhoods 
$\eta_\eps(E_D)$ and $\eta_\eps(\Gamma_D)$ 
as follows.  
Note that $E_D$ and $\Gamma_D$ can be viewed as subspaces of the products below
\begin{align*}
E_D \subset \LL_3^3 \x C_{q_1+q_2+q_3}\left[\coprod_1^3 \R\right] \x C_4[\R^3]  & & \Gamma_D \subset T_{\vec{n}} \LL_3^3 \x C_{q_1+q_2+q_3}\left[\coprod_1^3 \R\right] \x C_4[\R^3]
\end{align*}
where 
$q_1+q_2+q_3$ is determined by $D$.  
Let $\pi$ be the projection $C_4[\R^3] \to C_{q_1+q_2+q_3}[\R^3]$.  So $\pi$ is the identity if $D=L,M,$ or $R$, and $\pi$ forgets the ``free'' point if $D=T$.  
A point in $E_D$ is a point $(L, \mathbf{t}, \mathbf{x}) \in \LL_3^3 \x C_{q_1+q_2+q_3}[\coprod_1^3 \R] \x C_4[\R^3]$ with $L(\bt)=\pi(\bx)$.
Consider the subspace of such points with $|L(\bt) - \pi(\bx)| < \eps$, where the distance is measured using a metric on $C_{q_1+q_2+q_3}[\R^3]$.  Let $\eta_\eps(E_D)$ be the image of this subspace under the projection to $\LL_3^3 \x C_4[\R^3]$.

A point in $\Gamma_D$ is a point $((L_s)_{s\in \Delta^{n_1} \x ... \x \Delta^{n_3}}, \bt, \bx) \in T_{\vec{n}}\LL_3^3 \x C_{q_1+q_2+q_3}[\coprod_1^3 \R] \x C_4[\R^3]$ satisfying $L_{\widetilde{\gamma}(\bt)}(\bt) = \pi(\bx)$, where $\widetilde{\gamma}: C_{q_1+q_2+q_3}[\coprod_1^3 \R] \to \Delta^{n_1} \x ... \x \Delta^{n_3}$ is the map (\ref{GammaTilde}).
Consider the subspace of points $((L_s)_s, \bt, \bx)$ such that $| L_{\widetilde{\gamma}(\bt)}(\bt) - \pi(\bx) |<\eps$, where again the distance is measured in the compactified configuration space of $\R^3$.
Define $\eta_\eps(\Gamma_D)$ as the image of this subspace under the projection to $T_{\vec{n}}\LL_3^3 \x C_4[\R^3]$.

Clearly $\eta_\eps(E_D)$ and $\eta_\eps(\Gamma_D)$ are compatible in the sense that the middle vertical map in (\ref{EmbedIntoConfig}) takes the complement of $\eta_\eps(E_D)$ to the complement of $\eta_\eps(\Gamma_D)$.  Then because the two horizontal embeddings in the right-hand square of (\ref{EmbedIntoConfig}) are given by a product of the identity with the same embedding of $C_4[\R^3]$, we can find compatible neighborhoods of $E_D$ and $\Gamma_D$ in the right-hand column of (\ref{EmbedIntoConfig}).
For $D=T$, let $\widetilde{\eta}(E_T)$ and $\widetilde{\eta}(\Gamma_T)$ be such neighborhoods. 

For $D=L,M,R$, we can find compatible neighborhoods of $E_D \x S^2$ and $\Gamma_D \x S^2$ near the principal faces using those of $E_T$ and $\Gamma_T$ near their principal faces.  Away from these faces, the extensions of the embeddings to the products with $S^2$ are independent of the (punctured) links, so we can find compatible neighborhoods using those of $E_D$ and $\Gamma_D$.  Thus we obtain compatible neighborhoods $\widetilde{\eta}(E_D)$ and $\widetilde{\eta}(\Gamma_D)$ of $E_D \x S^2$ and $\Gamma_D \x S^2$ in the right-hand column of (\ref{EmbedIntoConfig}).
Appropriately gluing together the $\widetilde{\eta}(E_D)$ (respecively $\widetilde{\eta}(\Gamma_D)$) for the various $D$ gives compatible tubular neighborhoods of the horizontal embeddings in (\ref{EmbedIntoEuclid}).

Thus, for every sufficiently large $N$, we can collapse the complements of these neighborhoods. 
This yields the rows in the commutative square
\[
\xymatrix{
\LL_3^3 \x \left(\coprod_1^4 (\R^M \x [0,\infty)^N)/\sim\right) \ar[r]\ar[d] & E_g^{\nu_N} \ar[d] \\
T_n\LL_3^3 \x \left(\coprod_1^4 (\R^M \x [0,\infty)^N)/\sim\right) \ar[r] & \Gamma_D^{\nu_N}
}
\]
where the subscript $\nu_N$ denotes the Thom spaces of the normal bundles\footnote{Here and below, the symbol $\nu_N$ itself abusively stands for two different bundles, but we only use it together with the base space, so this should cause no confusion.  The same applies to the symbol $\nu$ in the next diagram.} of the horizontal embeddings in (\ref{EmbedIntoEuclid}).
This commutative diagram clearly descends to the quotients by boundaries.  Furthermore the cube of squares incorporating the indices $N$ and $N+1$ commutes.  Hence for $n\geq q$, there is a commutative square of spectra as below, on the left.  Applying cohomology (and suspension and Thom isomorphisms) and recalling that the fiber is 6-dimensional, we get the square on the right.
\begin{equation}
\label{SpectraSquare}
\xymatrix{
\Sigma^\infty (\LL_3^3)_+ \ar[r] \ar[d] & \E_g/ \d \E_g \ar[d]  & & H^{*-6}(\LL_3^3)  & H^*(E_g,\d E_g) \ar[l]  \\
\Sigma^\infty (T_{\vec{n}}\LL_3^3)_+ \ar[r] & \Gamma_g^\nu / \d \Gamma_g^\nu & & H^{*-6}(T_{\vec{n}}\LL_3^3) \ar[u] & H^* (\Gamma_g, \d \Gamma_g) \ar[l] \ar[u]  
}
\end{equation}

It remains to lift the class $[\beta] \in H^6(E_g, \d E_g)$ to $H^6(\Gamma_g, \d \Gamma_g)$.
Recall that $[\beta]$ was pulled back from $H^6(S^2 \x S^2 \x S^2, \mathcal{D})$ via a map $(E_g, \d E_g) \to (S^2 \x S^2 \x S^2, \mathcal{D})$, where $\mathcal{D}$ is a subspace containing the image of the faces at infinity $\d E_g$.  We first observe that there is a map $\Gamma_g \to S^2 \x S^2 \x S^2$.  In fact, each $\Gamma_D$ is a subspace of $T_{\vec{n}} \LL_3^3 \x C_4[\R^3]$, so we can use the same combination of ``unit vector'' and projection maps as for the $E_D(\x S^2)$ (see \cite[Section 5]{HoMilnor3ple}) to map $\Gamma_L \x S^2, \Gamma_M \x S^2, \Gamma_R \x S^2$, and $\Gamma_T$ to $S^2 \x S^2 \x S^2$.  Thus the diagram below commutes.
\begin{equation}
\label{MapToSpheres}
\xymatrix{
E_g \ar[r] \ar[d] & S^2 \x S^2 \x S^2 \\
\Gamma_g \ar[ur]  & 
}
\end{equation}
The boundary $\d \Gamma_g$ consists of configurations where some points have escaped to infinity.  The behavior at infinity for elements in $T_{\vec{n}}\LL_3^3$ is fixed in the same way as for links in $\LL_3^3$, so the image of $\d \Gamma_g$ in $S^2 \x S^2 \x S^2$ is the same as that of $E_g$.  (In other words, even though an element of $T_{\vec{n}}\LL_3^3$ is a family of punctured links, the variation of the links within this family is only inside the fixed compact subset of $\R^3$.)  Thus (\ref{MapToSpheres}) descends to a diagram of pairs, relative to the boundaries of $E_g$ and $\Gamma_g$ (faces at infinity) and the subspace $\mathcal{D} \subset S^2 \x S^2 \x S^2$ (which contains the image of these faces at infinity).  Applying cohomology and combining this with the right-hand square in (\ref{SpectraSquare}) yields the diagram in the Theorem statement and hence
completes the proof.
\end{proof}


As mentioned in the Introduction, we could also use integration of differential forms to deduce this result.  The advantage of our methods is the potential to yield integer or mod $p$ classes that cannot be realized via integration of forms.
The particular stage at which we see $\mu_{123}$ would be the same if we instead deduced this result using integration of differential forms.  In either case, we need to use four different configuration points, and on each strand, the maximum number of points over these four configuration spaces is 2.  Thus (2,2,2) seems to be the lowest stage at which one can obtain $\mu_{123}$ using configuration space integrals.    
While this matches Voli\'c's result that $T_n\K$ factors all rational type-$\frac{n}{2}$ knot invariants, Conjecture 1.1 of \cite{BCSS} is that $T_n \K$ factors all type-$(n-1)$ knot invariants.  (Our work in \cite{BCKS} provides some evidence for this conjecture.)
We might similarly expect to see the type-2 link invariant $\mu_{123}$ at a stage $(n_1, n_2, n_3)$ with $n_1+n_2+n_3=3$.  In fact, our work in \cite{CKKS} shows that link-homotopy classes of string links are distinguished at stage $(1,1,...,1)$ of the tower for \emph{link maps}. (This tower is defined similarly to $T_{\vec{n}}\LL_m^3$ but with the space of embeddings replaced by the space of link maps.)  
This implies that stage $(1,1,1)$ of that tower contains as much information as $\mu_{123}$ and the pairwise linking numbers.  A similar result is seen in the work of Munson \cite[Section 6.6]{MunsonMfdCalcLink}.  

\section{The aligned maps model}
\label{AlignedMaps}
In this section, we outline a possible alternative to the constructions above. 
This alternative approach uses a different model for the Taylor tower, namely a model of ``aligned maps,'' which is denoted $AM_n$ in the case of knots.  This is the main point:
\begin{quote}
The model $AM_n$ has a monoid structure compatible with stacking of long knots, while no such structure on the ``punctured knots model'' is known.
\end{quote}  
Thus $AM_n$ would be a more fruitful setting for studying in future work how our results interact with the monoid structure of stacking knots.  This alternative model still yields essentially the same main results via essentially the same proofs.  So we keep this Section brief, providing only an outline of the arguments.


\subsection{The monoid structure on the aligned maps model}
We showed in \cite{BCKS} that each stage of the Taylor tower for the space of long knots is an H-space with an operation compatible with connect-sum of long knots.  
There we used a variant of the aligned maps model where the $n$-th stage is roughly a space of maps of $\Delta^n \cong C_n \la I \ra$ to the simplicial compactification $C_n\la I^d\ra$ of configuration space.  (Here we are thinking of long knots as embeddings $I\incl I^d$ with fixed behavior at the endpoints; recall from Section \ref{Compactifications} that a configuration in $C_n \la I^d \ra$ has two fixed boundary points.)
However, for defining the monoid structure, we could have instead used the aligned maps model of \cite{SinhaTop, BCSS}, where the $n$-th stage $AM_n$ consists of maps from the associahedron $K_{n+2}\cong C_n[I]$ to the canonical compactification $C_n[I^d]$ of configuration space.  This latter option is needed for configuration space integrals because $C_n[M]$, unlike $C_n \la M \ra$, is a smooth manifold with corners.\footnote{Note that \cite{BCKS} and \cite{SinhaTop, BCSS} use the same notation $AM_n$ for these two slightly different variants.}  

In more detail, replacing strands of punctured knots by configurations requires us to record not only points in a configuration, but also a unit tangent vector at each point.  So let $C_n'[I^d]:=C_n[I^d] \x (S^{d-1})^n$, thought of as a subspace of $n+2$ points and unit tangent vectors, with two fixed boundary points $(\pm 1,0,...,0)$ and the fixed tangent vector $(1,0,...,0)$ at each of the two boundary points.
Then $AM_n$ is the space of maps $C_n[I] \to C_n[I^d]$ that are \emph{aligned} and \emph{stratum-preserving}, terms which we will now define.
Recall the stratifications on $C_n[I]$ and $C_n[I^d]$ by sets of subsets of $\{1,...,n\}$ from Section \ref{Compactifications}.  
Note that every set indexing a stratum in $C_n[I]$ also indexes a stratum in $C_n[I^d]$ (though not conversely, since the points in a configuration in $C_n[I]$ lie in order on the interval).
Then \emph{stratum-preserving} means that the image of a configuration in a given stratum in $C_n[I]$  must lie in the closure of the corresponding stratum in in $C_n[I^d]$.  The \emph{aligned} condition means that at a collision of points, the tangent vectors at all of the points involved must be equal.  
There is a canonical map $\K \to \Map(C_n[I], C_n'[I^d])$ given by evaluating a knot (and its derivative) on the $n$ points in a configuration in $C_n[I]$.  The resulting image is aligned and stratum-preserving, so this evaluation is a map $\K \to AM_n$.
See\cite[Definitions 5.3, 5.6; Lemma 5.7]{SinhaTop} or \cite[Definition 2.16]{BCSS} for more details.


\subsection{Bott--Taubes constructions on aligned maps models for knots and links}
To extend configuration space integral constructions to this aligned maps model $AM_n$, we make some slight adjustments to it.  First, we take only smooth maps, rather than continuous maps.  This allows the possibility of giving the mapping space a smooth structure (cf.~the beginning of Section \ref{BTIsOverTaylor}).  Also,
we need to replace $I$ by $\R$ to ensure the vanishing of integrals over faces at infinity (cf.~Remark \ref{KnotsInABox}):
\begin{definition}
\label{AMnRDefn}
Let $AM_n(\R)$ be the subspace of maps $\phi$ in $\Map(C_n[\R], C_n'[\R^d])$ which are aligned and stratum-preserving (as above) and also satisfy the following condition: 
\begin{quote}
if a point $x_i $ in a configuration $c$ is located outside of $(-1,1)$, then $p_i(\phi(c)) = ((x_i,0,...,0),(1,0,...,0))$, where $p_i : C_n'[\R^d] \to \R^d \x S^{d-1}$ is the projection to the $i$-th factor.
\end{quote}
\end{definition}
One can show that $AM_n(\R)$ is a model for $T_n \K$, just like $AM_n$ is, as in \cite[Section 6]{SinhaTop}.
One can also define an H-space structure as in \cite[Section 4]{BCKS} on $AM_n(\R)$.
We have
\[
AM_n(\R) \x C_n[\R] \subset 
\Map(C_n[\R], C_n'[\R^d]) \x C_n[\R] 
 \to C_n[\R^d]
\]
where the right-hand map above is the obvious canonical map followed by forgetting the tangent vectors. 
So we can write 
\[
\xymatrix{
\Gamma_{AM}[q;t] \ar[r] \ar[d] & C_{q+t}[\R^d] \ar[d] \\
AM_q(\R) \x C_q[\R] \ar[r] & C_q[\R^d] }
\]
where $\Gamma_{AM}[q;t]$ is defined as the pullback in this square.
We could also replace $AM_q(\R)$ by $AM_n(\R)$ for any $n\geq q$ using the projection map $AM_n(\R) \to AM_q(\R)$ in the tower.  Hence, there is a map $\Gamma_{AM}[q;t] \to AM_n(\R)$, and
as in Section \ref{BTonTaylorLinks}, we can show that this map is a bundle whose fibers are smooth manifolds with corners, isomorphic to those of the usual Bott--Taubes bundle over the space of long knots.
Thus we can integrate differential forms along the fiber.  This seems to give an alternative proof to Theorem \ref{FTKnotInvtsFactorThruTower}, thus bypassing the need for the somewhat technical construction of $\gamma$ in Voli\'{c}'s \cite[Proposition 4.1]{VolicConfig}, a map which we also used in Section \ref{BTonTaylorLinks}.

We can carry out neat embeddings and Pontrjagin--Thom collapse maps over $AM_n(\R)$, as we did in Section \ref{HoBTonTaylor} over the punctured knots model.
The collapse maps for $AM_n(\R)$ will be compatible (via the evaluation map) with the ones for the space of knots itself.  Thus the classes constructed in \cite{Rbo} factor through this model $AM_n(\R)$ for  sufficiently large $n$.  





Similar aligned maps models can be constructed in the case of long links.  
To model the stage $T_{\vec{n}} \LL^d_m = T_{(n_1,...,n_m)}\LL^d_m$, define $AM_{\vec{n}}\subset  \Map(C_{n_1}[\R] \x ...\x C_{n_m}[\R], C_{n_1+...+n_m}[\R^d])$ as the subspace of aligned, stratum-preserving smooth maps with prescribed behavior outside of $(-1,1)$, as in Definition \ref{AMnRDefn}.
Then for any $\vec{q}=(q_1,...,q_m)$ with $\vec{q} \leq \vec{n}$ (meaning every $q_i \leq n_i$), there is a projection $AM_{\vec{n}} \to AM_{\vec{q}}$.  So we can write the square 
\[
\xymatrix{
\Gamma_{AM}[q_1,...,q_m;t] \ar[r] \ar[d] & C_{q_1+...+q_m+t}[\R^d] \ar[d] \\
AM_{\vec{n}}(\R) \x C_{q_1+...+q_m}\left[\coprod_{i=1}^m\R\right] \ar[r] & C_{q_1+...+q_m}[\R^d] }
\]
where $\Gamma_{AM}[q_1+...+q_m;t]$ is defined as the pullback.  Thus we get bundles $\Gamma_{AM}[q_1,...,q_m;t] \to AM_{\vec{n}}(\R)$ for any $\vec{q} \leq \vec{n}$.  The fibers are isomorphic to those in the original 
Bott--Taubes bundle.  
Then the arguments in Sections \ref{BTonTaylorLinks} and \ref{HoMilnorOnTaylorSection} show that Theorems \ref{FTLinkInvtsFactor} and \ref{HoMilnorOnTaylor} also hold if one replaces the punctured knots model for $T_{\vec{n}} \LL^d_m$ by the aligned maps model $AM_{\vec{n}}(\R)$.


    \bibliographystyle{alpha}
    \bibliography{refs}

\end{document}